\documentclass[12pt]{amsart}
\usepackage{amscd,amssymb,graphics,lineno}
\usepackage{mathrsfs} 

\newtheorem{theorem}{Theorem}[section]
\newtheorem{corollary}[theorem]{Corollary} 
\newtheorem{lemma}[theorem]{Lemma}
\newtheorem{proposition}[theorem]{Proposition}
\theoremstyle{definition}
\newtheorem{definition}[theorem]{Definition}

\theoremstyle{remark}
\newtheorem{remark}[theorem]{Remark}

\newtheorem{example}[theorem]{Example}
\numberwithin{equation}{section}

\newcommand{\abs}[1]{\lvert#1\rvert}
\def\norm#1{\left\Vert#1\right\Vert}

\def\T {{\mathbb T}}
\def\N{{\mathbb N}}
\def\Z {{\mathbb Z}}

\def\R{{\mathbb R}}
\def\Ur{{\mathbb U}}
\def\Iso{{\mathrm{Iso}}}
\def\Emb{{\mathrm{Emb}}}
\def\Aff{{\mathrm{Iso}_{\mathrm{aff}}}}

\def\supp{{\mbox{\rm supp}\,}}
\def\e{\varepsilon}
\def\Homeo{{\mbox{\rm Homeo}\,}}
\def\Aut{{\mbox{\rm Aut}\,}}
\def\Lip{{\mbox{\rm Lip}\,}}
\def\orb{{\mbox{\rm orb}\,}}
\def\H {{\mathscr H}}

\newcounter{quest}
\stepcounter{quest}

\begin{document}

\title[Isometric actions on Banach spaces]
{Fixed point-free isometric actions of topological groups on Banach spaces}

\author[L. Nguyen Van Th\'e]{Lionel Nguyen Van Th\'e$^\dag$}

\address{Department of Mathematics and Statistics, University of Calgary, 2500 University Drive NW, Calgary, Alberta, Canada T2N 1N4}

\email{nguyen@math.ucalgary.ca}

\author[V. Pestov]{Vladimir G. Pestov$^\ddag$}

\address{Department of Mathematics and Statistics, 
University of Ottawa, 585 King Edward Ave., Ottawa, Ontario, Canada K1N 6N5}

\email{vpest283@uottawa.ca}

\thanks{$^\dag$The first named author gratefully acknowledges the support 
of the Department of Mathematics and Statistics Postdoctoral Program at the University of Calgary.}

\thanks{$^\ddag$The research of the second named author was partially supported by an NSERC discovery grant (2007---)
and University of Ottawa internal grants. The author thanks the Combinatorics group at the University of Calgary for hospitality extended in March 2008.}

\thanks{Both authors are grateful to the anonymous referee of this article for a number of useful suggestions.} 




\begin{abstract} 
We show that
every non-precompact topological group admits a fixed point-free continuous action by affine isometries on a suitable Banach space. Thus, precompact groups are defined by the fixed point property for affine isometric actions on Banach spaces. For separable topological groups, in the above statements it is enough to consider affine actions on one particular Banach space: the unique Banach space envelope $\langle\Ur\rangle$ of the universal Urysohn metric space $\Ur$, known as the Holmes space. At the same time, we show that Polish groups need not admit topologically proper (in particular, free) affine isometric actions on Banach spaces (nor even on complete metric spaces): this is the case for the unitary group $U(\ell^2)$ with strong operator topology, the infinite symmetric group $S_\infty$, etc.
\end{abstract}

\maketitle

\section{Introduction}

Kazhdan's property $(T)$ of locally compact, in particular discrete, groups
has found numerous uses, see \cite{dlHV} and especially \cite{BdlHV}. At the same time, the concept makes perfect sense for the most general Hausdorff topological groups, even if there are not many results known beyond the locally compact case, the interesting paper \cite{Bek3} being among a few exceptions.

Let $G$ be a topological group.
Recall that a strongly continuous unitary representation $\pi$ of $G$ in a Hilbert space $\H$ {\em has almost invariant vectors} if for every $\e>0$ and each compact $K\subseteq G$ there is a vector $\xi\in\H$ of norm one with the property $\norm{\xi-\pi_g\xi}<\e$ for all $g\in K$. Now $G$ has Kazhdan's property $(T)$ if every continuous unitary representation that has almost invariant vectors has an invariant vector of norm one. 

For second countable locally compact groups $G$ the property $(T)$ is equivalent to the following property known as $(FH)$: every continuous action of G by affine isometries on a Hilbert space has a fixed point. This result is known as {\em Delorme--Guichardet theorem,} see Ch. 2 in \cite{BdlHV}.
While the implication $(T)\Rightarrow (FH)$, due to Delorme \cite{delorme}, holds for every topological group, the converse implication $(FH)\Rightarrow (T)$, established by Guichardet \cite{guichardet} for sigma-compact locally compact groups, is in general invalid even for uncountable discrete groups. The following observation was made in print by Cornulier \cite{dC}, but noticed independently by some others. Let a group $G$ have Bergman's property \cite{bergman}, that is, whenever $G$ is represented as a union of an increasing countable chain $(W_i)_{i=1}^\infty$ of subsets, a suitable finite power of some $W_i$ equals $G$.
(Such is, for instance, the group $S_\infty$ of all self-bijections of natural numbers \cite{bergman}). Then $G$, viewed as a discrete group, has property $(FH)$. At the same time, every discrete group with property $(T)$ is finitely generated (\cite{BdlHV}, Th. 1.3.1.)

More generally, say, following Rosendal \cite{rosendal}, that a topological group $G$ has {\em property $(OB)$} if for every continuous action of $G$ on a metric space by isometries all orbits are bounded. ({\em Topological Bergman property} is another natural name for this property, suggested by the referee.) 
Such are discrete groups with Bergman's property, as well as topological groups $G$ bounded in the sense of Hejcman \cite{hejcman} and Atkin \cite{atkin} (that is, for every neighbourhood of identity $V$ in $G$, there are a finite subset $F\subset G$ and a natural $n$ with $FV^n=G$). The property $(OB)$ implies property $(FH)$ (\cite{rosendal}, Prop. 3.4), in direct consequence of the ``Lemma of the centre'' (\cite{BdlHV}, Section 2.2). However, for many known groups with property $(OB)$ it remains unknown whether or not they have property $(T)$. This is the case e.g. for the unitary group $U(\ell^2)$ with the uniform operator topology, the group $U(\ell^2(\Gamma))$ with the strong operator topology where $\abs{\Gamma}>\aleph_0$, the group $\Aut(X,\mu)$ of measure-preserving transformations of a standard Lebesgue measure space with the coarse topology, etc. At the same time, as shown recently by one of the present authors (V.P.) \cite{P08}, the properties $(T)$ and $(FH)$ are not equivalent even in the class of Polish (not necessarily locally compact) topological groups. For instance, the group $U_C(\ell^2)$ of all unitary operators which are compact perturbations of the identity, equipped with the uniform operator topology, has property $(OB)$ (hence $(FH)$), but is not a Kazhdan group ({\em ibid.})

Recently there has been some interest towards stronger versions of property $(FH)$ obtained if one allows continuous affine isometric actions of $G$ on Banach spaces from more general classes, especially $L^p$ spaces and other uniformly convex spaces \cite{BFGM,dCTV}.
In particular, some lattices in semisimple Lie groups have the fixed point property in this stronger sense, while others do not.

How far can one go in this direction by allowing affine actions of a topological group by isometries on more and more general classes of Banach spaces? Of course every affine action of a precompact group $K$ by isometries on a Banach space will always have a fixed point, obtained by integration on the completion of $K$, so the question really is: for what classes of Banach spaces do there exist {\em non-precompact groups} with the fixed point property for affine isometries?

Haagerup and Przybyszewska have shown \cite{HP} that every second countable locally compact non-compact group admits a metrically proper (in particular, fixed point-free) continuous affine action by isometries on a reflexive and strictly convex Banach space. Thus, the class of all reflexive Banach spaces is already too general: the corresponding notion for second countable locally compact groups is too strong and is never satisfied, save in the trivial case of compact groups. 

At the same time, one cannot hope to extend the result by Haagerup and Przybyszewska to non locally compact Polish groups. Indeed, by force of a theorem due to Megrelishvili \cite{megrelishvili01} stating that every weakly almost periodic function on the Polish group $\Homeo_+[0, 1]$ is constant, this particular group admits no nontrivial continuous affine actions by isometries on reflexive Banach spaces. (Cf. Corollary \ref{c:cannot} below.) The group $\Homeo_+[0, 1]$ therefore has the fixed point property for affine isometries on reflexive Banach spaces, but this in a sense happens for the wrong reason.

The main aim of this note is to ``dilute'' the property $(FH)$ inasmuch as possible, by allowing affine actions of {\em arbitrary} topological groups on {\em arbitrary} Banach spaces. In this context, we obtain an analogue of theorem of Haagerup and Przybyszewska, by observing that every topological group $G$ that is not precompact admits a continuous affine action by isometries on a Banach space without fixed points (Theorem \ref{th:1main}). In fact, this property characterizes precompactness. Notice that the action thus constructed can be chosen with bounded orbits (in contrast to actions on Hilbert spaces, for which a bounded orbit forces existence of a fixed point).

Moreover, there exists a single separable Banach space upon which every separable and non-precompact group $G$ acts continuously by affine isometries without fixed points. This is the remarkable {\em Holmes space} $\langle\Ur\rangle$ \cite{holmes}, the unique (up to an isometric isomorphism) Banach space spanned by the universal Urysohn metric space, $\Ur$ \cite{urysohn}, and containing it isometrically. In particular, among separable groups the precompact ones are characterized by the fixed point property with respect to affine isometric actions on the Holmes space $\langle\Ur\rangle$.

The proof uses a characterization of precompact groups obtained independently by Uspenskij (unpublished, cf. a remark in \cite{Usp98}) and Solecki \cite{solecki} as those topological groups G in which every neighbourhood of the identity, $U$, admits a finite set $F\subseteq G$ with $FUF=G$. Another component of the proof is the following observation of independent interest (Corollary \ref{c:extends}): every continuous action of a topological group $G$
by isometries on a metric space $X$ extends to an affine isometric action of $G$ on a suitable Banach space containing 
$X$ as a subspace and affinely spanned by it (the {\em Lipschitz-free} Banach space on $X$, cf. e.g. \cite{GC}). In order to extend the result to the Holmes space, we use the technique of Kat\v etov functions \cite{katetov}.

Simple cardinality arguments show that not every non-precompact topological group $G$ admits a free continuous action by affine isometries on a suitable Banach space. However, it turns out that even a Polish group need not admit a topologically proper (in particular, free) affine isometric action, moreover, this behaviour appears to be typical: examples of such groups include the infinite unitary group $U(\ell^2)$, the infinite symmetric group $S_\infty$, and a number of others (cf. the concluding Section \ref{s:non-existence}). In this sense, the full power of the result of Haagerup and Przybyszewska cannot be recovered for non locally compact Polish groups even if we allow actions on arbitrary Banach spaces.

\section{Preliminaries on actions of topological groups by isometries}

\subsection{Groups of isometries}
For a metric space, $\Iso(X)$ denotes the group of all surjective self-isometries of $X$ equipped with the topology of pointwise convergence, induced by the embedding $\Iso(X)\hookrightarrow X^X$ (or, which is the same,
compact-open topology). This group is second-countable if $X$ is such.

\subsection{Completeness}

The {\em left uniform structure} on a topological group $G$ is given by the basis of entourages of the diagonal consisting of all sets of the form
\begin{equation}
\label{eq:left}
V_L=\left\{(g,h)\colon g^{-1}h\in V\right\},\end{equation}
where $V$ runs over a neighbourhood basis of the identity element in $G$. Similarly, the {\em right uniform structure} is given by the basic entourages
\[V_R=\left\{(g,h)\colon gh^{-1}\in V\right\},\]
and the {\em two-sided uniform structure} is the supremum of the two, whose basis consists of the entourages
\[V_{\vee}=\{(g,h)\in G\colon g^{-1}h\in V\mbox{ and }gh^{-1}\in V\}.\]
A topological group $G$ is said to be {\em complete} if it is complete with regard to the two-sided uniformity. 

\begin{example} Let $X$ be a complete metric space. Then the topological group $\Iso(X)$ is complete.
\end{example}

Every topological group $G$ is isomorphic to an everywhere dense topological subgroup of a complete group $\hat G$, the completion of $G$.  It is easy to prove that every locally compact, in particular every compact, group is complete with regard to one-sided uniformities, but in general a topological group need not embed into a group complete with regard to the one-sided uniformities. This phenomenon, noticed in \cite{dieudonne}, will be discussed in some detail and used in Section \ref{s:non-existence}.

A general reference to uniformities on groups is \cite{RD}.

\subsection{Actions on metric spaces}
Let a topological group $G$ act by isometries on a metric space $X$. The following is an easy but useful observation.

\begin{lemma}
The following conditions are equivalent.
\begin{enumerate}
\item
The action $G\times X\to X$ is jointly continuous.
\item Every orbit map 
\[\orb_x\colon G\ni g\mapsto gx\in X,~~x\in X,\]
is continuous.
\item The orbit maps $\orb_x\colon G\to X$ are continuous for all $x$ from an everywhere dense subset of $X$.
\item The homomorphism $G\to \Iso(X)$ associated to the action of $G$ on $X$ is continuous.
\end{enumerate}
\qed
\label{l:orb}
\end{lemma}

\subsection{\label{ss:G/d}Transitive actions}
The following construction is a useful source of transitive isometric continuous actions of topological groups; moreover, every such action can be obtained in this way.

Let $d$ be a continuous left-invariant pseudometric on a topological group $G$. (According to the well-known Kakutani lemma, such pseudometrics determine the topology of $G$.) Set $H=\{g\in G\colon d(g,e)=0\}$. This $H$ is a closed subgroup of $G$. 
Denote by $X=G/H$ the quotient space equipped with the left action by $G$. 
The rule 
\[d_X(gH,kH) = d(g,k)\]
defines a metric on the quotient space $G/H$, because for every $g,k\in G$ and $h,h_1\in H$
\[d(gh,kh_1)\leq d(gh,g)+d(g,k)+d(k,kh_1)=d(g,k),\]
and similarly for the other inequality. 
Clearly, $d_X$ is translation invariant, and the action of $G$ on $X$ is continuous and transitive. We will denote $X=G/d$. 

Now suppose $G$ acts continuously and transitively by isometries on a metric space $X=(X,d_X)$. Select any point $\xi\in X$ and define for all $g,h\in G$
\[d_{\xi}(g,h)=d_X(g\xi,h\xi).\]
This $d_{\xi}$ is a continuous left-invariant pseudometric on $G$, and one can easily verify that the metric $G$-spaces $X$ and $G/d_{\xi}$ are isometrically isomorphic between themselves.

\subsection{Affine isometries}
The Mazur--Ulam theorem \cite{MU} says that every surjective isometry between two real Banach spaces is an affine map. In particular, every self-isometry of a Banach space $E$ is an affine isometry. To avoid confusion, we will not be using the symbol $\Iso(E)$ in case where $E$ is a Banach space, using $\Aff(E)$ instead for the group of all (affine) isometries of $E$. This way, continuous affine actions of $G$ on $E$ are in a one-to-one correspondence with continuous group homomorphisms $G\to\Aff(E)$. The symbol $\Iso_0(E)$ will stand for the group of all {\em linear} isometries of $E$ or, which is the same, all isometries of $E$ fixing the origin $0$.

The following is well-known (in the case of the Hilbert space, cf. Section 2.1 in \cite{BdlHV}).

\begin{proposition}
As a topological group, $\Aff(E)$ is isomorphic to the semidirect product $\Iso_0(E)\ltimes E$, where $E$ is viewed as an additive topological group with its norm topology, and the action of $\Iso_0(E)$ on $E$ is a tautological, or standard, one.
\end{proposition}

\begin{proof}
After the origin $0$ is chosen (and thus $E$ is made into a linear space),
an explicit isomorphism of topological groups $\Aff(E)\cong\Iso_0(E)\ltimes E$ is given by
\[\theta\colon \Aff(E)\ni g\mapsto (T_{-g(0)}\circ g, g(0))\in \Iso(E)\ltimes E,\]
where $T_a(x)=a+x$.
\end{proof}


\subsection{Precompact groups}
Recall that a topological group is {\em precompact} if it can be covered by finitely many left translates of every non-empty open subset. In other words, if $V$ is a neighbourhood of identity, there is a finite $F$ with $FV=G$. (For a treatment of uniform structures in topological groups, completions, precompactness etc., we recommend the book \cite{RD}.)

The following is essentially a well-known fact, although normally stated for compact groups only (cf. Lemma 2.3 in \cite{dCTV}).

\begin{proposition}
Let $G$ be a precompact topological group acting continuously by affine isometries on a Banach space $E$. Then $G$ admits a fixed point: for some $\xi\in E$ and all $g\in G$, one has
\[g\xi=\xi.\]
\label{p:fpbs}
\end{proposition}

\begin{proof}
By an equivalent definition, a topological group $G$ is precompact if and only if its completion $\hat G$ is compact. The continuous action $G\to\Aff(E)$ extends to a continuous homomorphism from $\hat G$ to $\Aff(E)$ (a complete group), that is, a continuous action by affine isometries of the compact group $\hat G$. 
Choose any point $\xi\in E$. The orbit $\hat G\xi$ is a compact metric space, and so is the convex hull $C={\mathrm{conv}}\,\hat G\xi$. (The convex hull only depends on the affine structure of $E$, so one can turn $E$ into a Banach space in an arbitrary way and use a well-known result valid for all complete locally convex spaces,
see e.g. \cite{schaefer}, Corollary on p.50.)
Moreover, $C$ is invariant under the action of $\hat G$ by affine maps.
But every compact group admits a fixed point in every convex compact set upon which it acts by affine homeomorphisms (this is a reformulation of a statement that every compact group is amenable).
\end{proof}

\section{Affine isometric actions on reflexive Banach spaces}

Recall that a bounded continuous scalar-valued function $f$ on a topological group $G$ is {\it weakly almost periodic} \cite{Ru} if the orbit of $f$ with regard to the left regular representation of $G$ is weakly relatively compact in the
Banach space $CB(G)$ of all bounded continuous functions, equipped with the supremum norm. 

Weakly almost periodic functions on a topological group $G$ are closely related to strongly continuous representations of $G$ by isometries in reflexive Banach spaces. Let $\pi\colon G\to\Iso_0 (E)$ be such a representation. Fix an element $\xi\in E$ and a bounded linear functional $\phi\in E^\prime$. Then the function
\[G\ni g\mapsto \phi(\pi_g\xi)\]
(a {\em matrix coefficient} of $\pi$)
is weakly almost periodic. This follows easily from the fact that the unit ball of a reflexive Banach space is weakly compact. 
A considerably less trivial observation is that {\em all} weakly almost periodic functions on $G$ arise this way \cite{shtern,megrelishvili99}. 

Every locally compact group $G$ admits a rich collection of weakly almost periodic functions: they separate points and closed subsets of $G$. (Indeed, this is true even of elements of the Fourier-Stiltjes algebra on $G$, that is, coefficients of unitary representations, cf. e.g. \cite{gao}.)  
However, for more general Polish groups this is no longer the case. Denote by $\Homeo_+[0,1]$ the group of endpoint-preserving homeomorphisms of the closed unit interval, equipped with the compact-open topology. This is one of the best studied examples of Polish groups. As shown by Megrelishvili \cite{megrelishvili01}, this group admits no non-constant continuous weakly almost periodic functions. (Another such example \cite{P07}, which is really a consequence of the above one, is the group $\Iso(\Ur_1)$ of isometries of a sphere in the universal Urysohn metric space $\Ur$, an object which will be defined later in our paper.)

Thus, there exist Polish groups (such as $\Homeo_+[0,1]$ or $\Iso(\Ur_1)$) that admit no nontrivial strongly continuous representations by isometries in reflexive Banach spaces \cite{megrelishvili01}.
We are going to strengthen this observation a little bit.

\begin{proposition}
Suppose a topological group $G$ admits no non-constant continuous weakly almost periodic functions. Then the only continuous action of $G$ by affine isometries on a reflexive Banach space $E$ is the trivial (constant) action.
\end{proposition}

\begin{proof} Let $\pi\colon G\to\Aff(E)$ be a continuous homomorphism corresponding to an affine isometric action.
Denote by $q\colon \Aff(E)\cong \Iso_0(E)\ltimes E\to\Iso_0(E)$ the quotient map (the first coordinate projection). The composition $q\circ\pi$ is a continuous representation of $G$ by linear isometries of $E$ and so is trivial, meaning that $\pi$ maps $G$ into the kernel of $q$, that is, $\pi(G)\subseteq E$. However, continuous one-dimensional unitary representations separate points of $E$, as is seen by associating to every $\phi\in E^\prime$ a character $e^{2\pi\mathrm{i}\phi}$. This allows to conclude: $\pi(G)$ is the identity.
\end{proof}

\begin{corollary}
The Polish topological groups $\Homeo_+[0,1]$ and $\Iso(\Ur_1)$ cannot act by affine isometries on a reflexive Banach space in a non-trivial way.
\label{c:cannot}
\end{corollary}

In particular, a result of Haagerup and Przybyszewska \cite{HP} cited in the Introduction does not extend to non locally compact, non precompact Polish groups. 

\section{Fixed point-free actions on metric spaces}

\begin{theorem}[Uspenskij (unpublished, cf. \cite{Usp98}); Solecki \cite{solecki}, Lemma 1.2]
A topological group $G$ is precompact if and only if for every non-empty open subset $V$ of $G$ there is a finite set $F$ with the property 
\[FVF=G.\]
\qed
\end{theorem}

Notice that the condition is formally weaker than the usual definition of precompactness. For a simple proof of the above criterion, see \cite{BT}, Proposition 4.3.

\begin{definition}
Say that an action of a group $G$ by isometries on a metric space $X$ is {\em strongly moving} {\em with constant} $\e_0 > 0$ if for every finite subset $F\subseteq X$ there is a $g\in G$ with the property
\begin{equation}
\label{eq:su}
d_X(F,gF)\geq \e_0.\end{equation}
(Here the distance is understood as a distance between two sets.)
\end{definition}

\begin{remark}
In the above definition, one can replace ``every finite subset $F\subseteq G$'' with ``every compact subset $K\subseteq G$'' by approximating $K$ with finite nets.
\end{remark}

\begin{corollary} For a topological group $G$, the following conditions are equivalent.

\begin{enumerate}
\item $G$ is non-precompact.
\label{1}
\item $G$ admits a continuous transitive action by isometries on a non empty metric space $X$ that is strongly moving.
\label{2}
\item $G$ admits a continuous strongly moving action by isometries on a non empty metric space $X$. 
\label{3}
\end{enumerate}

\label{c:precomp}
\end{corollary}

\begin{proof}
(\ref{1})$\Rightarrow$(\ref{2}).
Suppose that $G$ is not precompact. According to the Solecki--Uspenskij criterion, there exists an open neighbourhood $V$ of the identity with the property that for every finite set $F\subseteq G$, one has $FVF\neq G$. Choose a continuous left-invariant pseudometric $d$ on $G$ whose open unit ball around identity, $O^d_1(e)$, is contained in $V$. Form a metric space $X=G/d=G/H$ equipped with an action of $G$ by isometries as in Paragraph \ref{ss:G/d}.

Given a finite set $F\subseteq X$, choose a finite symmetric $\Phi\subseteq G$ so that the image of $\Phi$ under the quotient map $G\to G/H$ contains $F$. Choose
$g\in G\setminus \Phi V\Phi$.
We claim that Eq. (\ref{eq:su}) holds with $\e_0=1$. Indeed, assuming that for some $x,y\in \Phi$ one has $d_X(xH,gyH)<1$, one concludes: $d(x^{-1}gy,e)=d(x,gy)<1$ and so
$x^{-1}gy\in V$, that is, $g\in xVy^{-1}\subseteq \Phi V\Phi$, a contradiction.

(\ref{2})$\Rightarrow$(\ref{3}). Trivial.

(\ref{3})$\Rightarrow$(\ref{1}).
If $G$ is precompact, then so is every orbit of $X$. Denoting for every $\e>0$ by $F_\e$ a finite $\e$-net in any chosen $G$-orbit, one has $d_X(F_\e,gF_\e)<\e$ for all $g\in G$, meaning the action is not strongly moving.
\end{proof}

\section{Lipschitz-free Banach spaces}

It is well known that every metric space embeds isometrically into a normed space. Moreover, such an embedding can be performed in a universal way. We want to revise details of this construction.

Let $X=(X,d)$ be a metric space and $\ast$ a point in $X$. We will refer to the triple $(X,d,\ast)$ as a {\em pointed metric space.} 
Denote by $L(X,\ast)$, or simply
by $L(X)$, a real vector space having $X\setminus\{\ast\}$
as its Hamel basis and $\ast$ as zero. 
Equip $L(X)$ with the largest prenorm $p$ whose restriction to $X$ is bounded by the distance $d$: for all $x,y\in X$, one has $p(x-y)\leq d(x,y)$. 
Then it is easy to show that $p$ is a norm, inducing the distance $d$ on $X$ \cite{K-R,AE,michael,Fl1,Fl2,GC}. 

Denote by $\Lip(X,\ast)$ the linear space of all Lipschitz
functions $f\colon X\to\R$ with the property $f(\ast)=0$, equipped with the norm
$\norm f$ equal to the infimum of all Lipschitz constants for $f$. 
For an $x\in X$, denote by $\hat x$ the usual evaluation
functional $\Lip(X,\ast)\ni f\mapsto f(x)\in \R$. The mapping 
\[X\ni x\mapsto \hat x\in \Lip(X,\ast)^\prime\]
is an isometric embedding of $X$ into the dual Banach space of
$\Lip(X,\ast)$, and the image of $X\setminus\{\ast\}$ is linearly independent.

To every element of $L(X,\ast)$, viewed as a finitely-supported measure on
$X\setminus\{\ast\}$, one can associate a bounded linear functional
on $\Lip(X,\ast)$, and this determines an embedding of $L(X,\ast)$ into the dual space $\Lip(X,\ast)^\prime$ as a normed subspace. It is an easy exercise to show that the norm induced on $L(X)$ from $\Lip(X,\ast)^\prime$ is the maximal
prenorm that we are after. Notice also that $X$ is closed in
$L(X)$. 

Another way to define $L(X)$ is through the following universal property. 

\begin{theorem}
\label{univer} Let $(X,\ast)$ be a pointed metric space,
let $E$ be a normed space, and let $f\colon X\to E$ be a
1-Lipschitz map with the property $f(\ast)=0$. Then there is
a unique linear operator $\bar f\colon L(X)\to E$ of norm $1$
extending $f$.  
\end{theorem}

\begin{proof}
There is a unique linear mapping $\bar f\colon L(X)\to E$ extending $f$, and it only remains to notice that the prenorm $N(x) = \norm{\bar f(x)}_E$ has the property that its restriction to $X$ is bounded by $d$.
\end{proof}

The Banach space completion of $L(X)$ is
denoted by $B(X)$ and called the ({\em Lipschitz})
{\it free Banach space} on $(X, d, \ast)$. It has a universal property
similar to Theorem \ref{univer}, but with respect to
all {\it Banach} spaces $E$.

For instance, in the case where $X=\Gamma\cup\{\ast\}$ is a set equipped with a  
$\{0,1\}$-valued metric, 
the free Banach space $B(\Gamma\cup\{\ast\})$ is just $\ell^1(\Gamma)$.

For various choices of distinguished point $\ast\in X$, the resulting free normed (Banach) spaces $L(X,\ast)$ (resp., $B(X,\ast)$) are isometrically isomorphic between themselves. This is seen by applying the universality property to the 1-Lipschitz mapping 
\[B(X,\ast)\supset X\ni x\mapsto x-\ast+\star\in B(X,\star).\]

All of the above facts are very well-known, and in addition to the references given above a useful summary can be found in \cite{GP}.

Here we observe that a similar universal property holds with regard to all isometries (not necessarily preserving the distinguished point) and affine isometric mappings. 

\begin{proposition}
Let $(X,\ast)$ be a pointed metric space, let $E$ be a normed space, and let $f\colon X\to E$ be a $1$-Lipschitz mapping. There exists a unique affine $1$-Lipschitz mapping $\tilde f\colon L(X,\ast)\to E$ whose restriction to $X$ is $f$.
\label{th:affine}
\end{proposition}

\begin{proof}
The mapping $g(x) = f(x) - f(\ast)$ is a $1$-Lipschitz mapping from $X$ into $E$ sending the distinguished point $\ast$ to zero, and so $g$ extends to a $1$-Lipschitz linear mapping $\bar g\colon L(X,\ast)\to E$. For all $x\in L(X,\ast)$, define $\tilde f(x) = \bar g(x) + f(\ast)$. If $x\in X$, one has $\tilde f(x) = f(x) - f(\ast) +f(\ast) = f(x)$, and this $\tilde f$ is affine and $1$-Lipschitz. If $h$ is another affine $1$-Lipschitz mapping from $L(X,\ast)$ to $E$ with $h\vert_X=f$, the formula $\hat h(x) = h(x) - h(\ast)$ defines a linear $1$-Lipschitz map whose restriction to $X$ equals $g$, meaning that $h(x)=\tilde f(x)$.
\end{proof}

By proceeding to completions, we obtain:

\begin{corollary}
Let $(X,\ast)$ be a pointed metric space, let $E$ be a Banach space, and let $f\colon X\to E$ be a $1$-Lipschitz mapping. There exists a unique affine $1$-Lipschitz mapping $\tilde f\colon B(X,\ast)\to E$ whose restriction to $X$ is $f$. \qed
\end{corollary}

\begin{corollary}
Let $(X,\ast)$ be a pointed metric space. Then every surjective self-isometry of $X$ extends in a unique way to a surjective affine self-isometry of $L(X,\ast)$ (and of $B(X,\ast)$).
\label{c:extends}
\end{corollary}

\begin{proof}
Both $f$ and $f^{-1}$ extend in a unique fashion to affine $1$-Lipschitz self-maps of $L(X)$ (respectively, of $B(X)$). Denote them by $\tilde f$ and $\hat f$. Because of uniqueness, these two maps should be mutually inverse and so are surjective affine isometries.
\end{proof}

Now let $G$ be a group acting by isometries on a metric space $X$. Denote the action by $\tau$. Choose a distinguished point $\ast\in X$.
For every $g\in G$, the self-isometry $\tau_g$ of $X$ extends to a unique affine self-isometry $\tilde\tau_g$ of the free Banach space $B(X,\ast)$. This defines an action of $G$ by affine isometries on $B(X,\ast)$, because of the same uniqueness consideration: $\widetilde{\tau_{fg}}$ must coincide with $\tilde\tau_f\circ \tilde\tau_g$.

\begin{proposition}
Let $G$ be a topological group acting on a metric space $X$ continuously by isometries. Then the affine extension of the action over $B(X,\ast)$ is continuous as well.
\label{p:ext}
\end{proposition}

\begin{proof}
According to Lemma \ref{l:orb}, it is enough to verify that every orbit map 
\[\orb_x\colon G\ni g\mapsto \tilde \tau_g(x)\in B(X,\ast)\]
is continuous, where $x\in L(X,\ast)$. This in its turn follows easily from the (presumed) continuity in the case where $x\in X$ because $B(X)$ is the affine subspace spanned by $X$.

\end{proof}

\begin{remark}
To see why the above result is interesting, notice that in the linear case an action of a group $G$ on a metric space by isometries can be extended to an action of $G$ on $B(X)$ the way it was done in \cite{P86} if and only if $G$ possesses bounded orbits in $X$, cf. \cite{schroeder}. No such restriction is necessary in the affine case.
\end{remark}

\begin{theorem}
Let a group $G$ act on a metric space $X$ by isometries. Suppose that the action is strongly moving. Then the extension of the action to an affine action on $B(X,\ast)$ is strongly moving with the same constant.
\label{th:extension}
\end{theorem}

\begin{proof} It is enough to verify that the action of $G$ on $L(X,\ast)$ is strongly moving with the same constant $\e_0$ as the constant of the action on $X$, because $L(X,\ast)$ is dense in $B(X,\ast)$. 
Let $F\subseteq L(X,\ast)$ be a finite subset of vectors. Denote 
\[\Phi = \{\ast\}\cup\bigcup_{v\in F}\supp v\]
and choose a $g\in G$ with the property
\[d_X(\Phi,g\Phi)\geq\e_0.\]
Let $v,w\in F$. Write
\[v=\sum_i\lambda_i x_i,~~w=\sum_i\mu_ix_i,\]
where $x_i\in\Phi$ and the sums are finite. According to the way the affine mapping $\tilde g$, extending $g$, was defined in Proposition \ref{th:affine},
\begin{eqnarray*}
\tilde g(v) &=& \sum_i\lambda_i (g(x_i) - g(\ast))+ g(\ast) \\
&=& \sum_i\lambda_i g(x_i)+ \left( 1- \sum_i\lambda_i\right) g(\ast).\end{eqnarray*}
There is a $1$-Lipschitz real function $h$ on $X$ taking value $0$ at $\Phi$ (in particular, vanishing at $\ast$) and the value $\e_0$ on $g\Phi$, for example,
\[h(x) = \min\left\{\e_0,d(x,\Phi)\right\}.\]
Looking at $w$ and $\tilde g v$ as linear functionals on the space $\Lip(X,\ast)$, one observes that 
\[w(h) = \sum_i\mu_i h(x_i) =0,\]
while
\[\tilde g(v)(h) = \sum_i\lambda_i h(gx_i) + \left( 1- \sum_i\lambda_i\right)h(g(\ast)) =
\e_0,
\]
meaning that $\norm{w-\tilde gv}\geq\e_0$.
\end{proof}

By combining Theorem \ref{th:extension} with Corollary \ref{c:precomp} and Proposition \ref{p:fpbs}, we obtain the first main result of this article:

\begin{theorem}
For a topological group $G$, the following conditions are equivalent.
\begin{enumerate}
\item $G$ is non-precompact.
\item $G$ admits a fixed point-free continuous action by affine isometries on a Banach space.
\item $G$ admits a strongly moving continuous action by affine isometries on some Banach space.
\end{enumerate}
\label{th:1main}
\qed
\end{theorem}

\begin{remark}
The action as above can be always chosen with bounded orbits (in contrast to actions on Hilbert spaces, for which a bounded orbit forces existence of a fixed point). This follows from the proof of Corollary \ref{c:precomp}, because $d$ can be chosen bounded. Moreover, one cannot expect in general an unbounded action, as the non-precompact group may satisfy Bergman's property, e.g. $S_\infty$ with its standard Polish topology.
\end{remark}
%


\section{Urysohn space and Holmes space}

The purpose of this section is to prove that when $G$ is separable (or, more generally, $\omega$-bounded in the sense of \cite{guran}), Theorem \ref{th:1main} can actually be strengthened by considering actions of $G$ on one single particular Banach space. Using the terminology of the previous section, this space is obtained by taking the free Banach space over a particular complete separable metric space known as the Urysohn space. 

The Urysohn space was constructed as a response to the question posed by Fr\'echet: is there a separable metric space $X$ universal for the class of all separable metric spaces (that is, into which any separable metric space embeds isometrically)? The answer came in 1925 with the construction by Urysohn \cite{urysohn} of a complete separable metric space $\Ur$ with the required property, as well as the following {\em ultrahomogeneity property}: every isometry $f$ between finite isometric subspaces $X$ and $Y$ of $\Ur$ can be extended to an isometry $\hat{f}$ of $\Ur$ onto itself. Moreover, Urysohn also showed that $\Ur$ is the unique complete separable metric space that is both universal and ultrahomogeneous. 
For an introduction to the theory of the Urysohn space, see e.g. \cite{melleray}, or Chapter 5 in \cite{P06}, or \cite{NVT}.

Much later, when studying how $\Ur$ may be embedded into the Banach space $C[0,1]$ of all continuous functions from $[0,1]$ to $\R$ equipped with the sup norm \cite{holmes}, Holmes discovered that $\Ur$ has the following remarkable property: for every isometric embedding $i$ (resp. $j$) of $\Ur$ into a Banach space $Y$ (resp. $Z$) such that the element $0_X$ (resp. $0_Y$) is in the range of $i$ (resp. $j$), the closed linear span of $i(\Ur)$ in $Y$ is linearly isometric to the closed linear span of $j(\Ur)$ in $Z$. In other words, up to linear isometry, the space $\Ur$ generates a unique Banach space, today called  {\em Holmes space} and denoted $\langle\Ur\rangle$. Now, the free Banach space $B(\Ur)$ also being generated by $\Ur$, the space $B(\Ur)$ actually {\em is} $\langle\Ur\rangle$. Note that several very elementary questions are still open concerning $\langle\Ur\rangle$. For example, does it have a basis? On the other hand, what is now well understood, is the property that allows $\Ur$ to generate a unique Banach space (see the recent paper by Melleray, Petrov and Vershik \cite{mvp}, or else \cite{P06}, p. 113).

For our purposes, the space $\langle\Ur\rangle$ is relevant because of the following specialized version of Theorem \ref{th:1main}. 

Recall that a topological group $G$ is {\em $\omega$-bounded} if it can be covered by countably many translates of every non-empty open subset. This condition is equivalent to the following: $G$ embeds as a topological subgroup into the direct product of a family of second-countable groups with the standard product topology. See \cite{guran}.

\begin{theorem}
For a separable (more generally, $\omega$-bounded) topological group $G$, the following conditions are equivalent.
\begin{enumerate}
\item $G$ is non-precompact.
\item $G$ admits a fixed point-free continuous action by affine isometries on the Holmes space.
\item $G$ admits a strongly moving continuous action by isometries on the Holmes space.
\end{enumerate}
\label{th:2main}
\qed
\end{theorem}

\begin{corollary}
For a separable (more generally, $\omega$-bounded) topological group $G$, the following conditions are equivalent:
\begin{enumerate}
\item $G$ is precompact.
\item Every continuous affine action of $G$ on the Holmes space $\langle\Ur\rangle$ by isometries has a fixed point.
\end{enumerate}
\qed
\end{corollary}

The implication $(3)\Rightarrow (2)$ in Theorem \ref{th:2main} is trivial, while $(2) \Rightarrow (1)$ comes from the easy (sufficiency) part of Corollary \ref{c:precomp}. In order to establish $(1)\Rightarrow (3)$, we are going to prove that if $G$ is not precompact, then it admits a strongly moving action on $\Ur$, and hence on $\langle\Ur\rangle = B(\Ur)$ by virtue of Theorem \ref{th:extension}. 

The required action on $\Ur$ is obtained thanks to the construction of $\Ur$ by Kat\v{e}tov in \cite{katetov}. The main concept of this construction is the concept of {\em Kat\v{e}tov function} (which
functions were introduced and used earlier by Flood in \cite{Fl1}, \cite{Fl2}, although in a somewhat different context). 

\begin{definition}
Let $X = (X, d)$ be a metric space and $f \colon X \to [0,+\infty)$. Say that $f$ is {\em Kat\v{e}tov over $X$} when for every $x, y \in X$, \[ |f(x) - f(y)| \leq d(x,y) \leq f(x) + f(y).\] 
\end{definition}

Equivalently, $f$ is Kat\v{e}tov over $X$ when it can be used to define a 1-point metric extension $X\cup \{p\}$ of $X$
by setting for every $x$ in X, $d (x, p) = f(x)$. In the sequel, the set of all Kat\v{e}tov maps over $X$ is written $E(X)$. 

With respect to the Urysohn space, the crucial fact is that those functions provide another characterization of $\Ur$: up to isometry, the space $\Ur$ is the unique complete separable metric space such that for every finite $F \subseteq \Ur$ and $f \in E(F)$, there is $y \in \Ur$ such that for every $x \in F$, $d(x,y)=f(x)$.

The following method (due to Kat\v{e}tov) allows to construct a space where that condition holds. Starting from any separable metric space $X$, construct first a metric space $X^*$ containing $X$ as a subspace (or, more exactly, into which $X$ embeds isometrically) and where for every finite $F \subseteq X$ and every $f \in E(X)$, there is $y \in X^*$ such that for every $x \in F$, $d(x,y)=f(x)$. Next, set $X_0 = X$ and define inductively $X_{n+1} = X_n ^*$. The completion $\overline{\bigcup _{n=0} ^{\infty} X_n}$ of the union $\bigcup _{n=0} ^{\infty} X_n$ is then separable and satisfies the required condition. 

As for $X^*$, a possible way to construct it from $X$ is as follows: for $Y \subseteq X$ finite and $f \in E(Y)$, $f$ can be used to define a $1$-point metric extension $\hat{f}$ of $X$ by setting \[ \hat{f}(x) = \min_{y\in Y} f(y) + d(y,x) .\] 

Now, set \[ X^* = \bigcup \{\hat{f}: Y\subseteq X \ \textrm{finite}, \ f \in E(Y) \}.\]

Then $X^*$ equipped with the sup metric is as required: first, it embeds $X$ isometrically because every $x$ may be identified with the function $\hat{f_x}$ where $f_x$ is the unique Kat\v{e}tov function on $\{ x\}$ with value $0$. Next, given any finite $F \subseteq X^*$ and $f \in E(F)$, there is indeed $y \in X^*$ such that $d(x,y) = f(x)$ for every $x \in F$: simply take $y = \hat{f}$.     

Kat\v{e}tov method is not only relevant for the previous construction of $\Ur$. As noted by Uspenskij in \cite{Usp98}, it also allows to prove that any separable metric space $X$ embeds into $\Ur$ in such a way that every $g \in \Iso (X)$ extends to $\hat{g} \in \Iso (\Ur)$ with the map $\Iso(X) \ni g \mapsto \hat{g} \in \Iso (\Ur)$ being a topological group embedding. Indeed, consider the chain \[ X \subseteq X^* = X_1 \subseteq X_1 ^* = X_2 \subseteq \ldots \subseteq \Ur.\]

Notice that the inclusion $X \subseteq X^*$ is such that every $g \in \Iso (X)$ extends to an embedding $g^* \in \Iso (X^*)$ in such a way that $\Iso(X) \ni g \mapsto g^* \in \Iso (X^*)$ is a topological group embedding: if $g \in \Iso (X)$ and $f \in X^*$, set $g^* (f) = f\circ g^{-1}$. Iterating this process, $\Iso(X)$ embeds into $\Iso(\overline{\bigcup _{n=0} ^{\infty} X_n}) = \Iso(\Ur)$. For $g \in \Iso(\Ur)$, the image via this embedding is obtained by extending $g$ to an element $g_n$ of $\Iso (X_n)$ so that $g_{n+1} \vert _{X_n} = g_n$ for every $n$ and by considering the isometry $\hat{g}$ induced on $\overline{\bigcup _{n=0} ^{\infty} X_n} = \Ur$ by the sequence $(g_n)_n$.

For our purposes, Kat\v{e}tov's construction is important because given a topological group $G$, it allows to make a strongly moving action of $G$ on $\Ur$ out of any strongly moving action of $G$ by isometries on an arbitrary separable metric space $X$. This is how we proceed: 

\begin{proposition}
Let $Y$ be a metric space, $A, B \subseteq Y$ finite with $d(A,B) \geq \varepsilon$, $\varphi \in E(A)$ and $\psi \in E(B)$. Then $\|\hat{\varphi} - \hat{\psi} \| \geq \varepsilon$. 
\label{prop:K}
\end{proposition}   

\begin{proof}
If $| \hat{\varphi} - \hat{\psi} |$ takes a value larger than $\varepsilon$ on $A$, we are done. Otherwise, for every $a \in A$ \[ | \hat{\varphi}(a) - \hat{\psi}(a)| < \varepsilon\] and so \[ \hat{\psi}(a) - \varepsilon < \hat{\varphi}(a).\] But by definition of $\hat{\psi}$, \[ \hat{\psi}(a) = \min_{b \in B} \psi(b) + d(b,a) \geq \min_{B} \psi + \varepsilon.\] So for every $a\in A$ \[\min_{B} \psi < \hat{\varphi}(a) (= \varphi(a)).\] Choose now $b_0 \in B$ where $\psi$ reaches its minimum. Then \[ \hat{\varphi}(b_0) = \min_{a\in A} \varphi(a) + d(a,b_0) > \min_{B} \psi + \varepsilon = \hat{\psi}(b_0)+\varepsilon.\] So  \[  |\hat{\varphi}(b_0) - \hat{\psi}(b_0)|>\varepsilon. \]
\end{proof}

We saw above that we could embed $\Iso(X)$ into $\Iso(X^*)$ by setting, for $g \in \Iso (X)$ and $\hat{f} \in X^*$, $g^* (\hat{f}) = \hat{f}\circ g^{-1}$. Therefore, if a topological group $G$ acts on a metric space $X$, then we can extend the action of $G$ to an isometric action on $X^*$ by setting, for $g \in G$ and $\hat{f} \in X^*$, $g\hat{f}=\hat{f}\circ g^{-1}$.  

\begin{proposition}
Let $G$ be a topological group acting on a metric space $X$. Suppose that the action is strongly moving with constant $\varepsilon$. Then the extension of the action of $G$ to an isometric action on $X^*$ is strongly moving with constant $\varepsilon$.  
\end{proposition}

\begin{proof}
Fix $H \subseteq X^*$ finite. We need to find $g\in G$ such that $d(H,g H) \geq \varepsilon$. To achieve that, fix $F \subseteq X$ finite but large enough such that every element of $H$ is of the form $\hat{f}$ for some $f \in E(F)$. Observe that for every $f \in E(F)$, $gf \in E(gF)$ and that $g\hat{f} = \widehat{gf}$. Indeed, for $x \in X$,  

\begin{eqnarray*}
g\hat{f}(x)&=& \hat{f}(g^{-1}x)\\
&=& \min_{y \in F} f(y) + d(y,g^{-1}x)\\
&=& \min_{y\in gF} f(g^{-1}y) + d(g^{-1}y,g^{-1}x)\\
&=& \min_{y\in gF} gf(y) + d(y,x)\\
&=& \widehat{gf}(x).  
\end{eqnarray*}

It follows that every element of $gH$ is of the form $\hat{f}$ for some $f \in E(gF)$, while every element of $H$ is of the form $\hat{f}$ for some $f \in E(F)$. Therefore, in virtue of Proposition \ref{prop:K}, it is enough to choose $g\in G$ such that $d(F,gF)\geq \varepsilon$. 

\end{proof}

Iterating the previous proposition, we obtain:  

\begin{proposition}
Let $G$ be a topological group acting on a metric space $X$. Suppose that the action is strongly moving with constant $\varepsilon$. Then the extension of the action of $G$ to an isometric action on $X_n$ is strongly moving with constant $\varepsilon$.  
\end{proposition} 

In the case where $X$ is separable, the fact that $\Ur$ may be seen as the completion of $\bigcup_{n=0}^{\infty} X_n$ consequently allows to prove: 

\begin{proposition}
Let $G$ be a topological group acting on a metric space $X$. Suppose that the action is strongly moving with constant $\varepsilon$. Suppose also that $X$ is separable. Then the extension of the action of $G$ to an isometric action on $\Ur$ is strongly moving with constant $\varepsilon$. 
\end{proposition}

\begin{proof}
Fix $H \subseteq \Ur$ finite. Then there is $n \in \N$ such that $H \subseteq X_n$ and by the previous proposition, we can find $g \in G$ such that $d(H,gH)\geq \varepsilon$.
\end{proof}

To finish the proof of Theorem \ref{th:2main}, it remains to show that if $G$ is separable and not precompact, then it admits a strongly moving action by isometries on a separable metric space. This is done by observing that the space $X$  provided by Corollary \ref{c:precomp} is separable in the case where $G$ is.
Note that if more generally, $G$ is $\omega$-bounded, then the space $X$ constructed in Corollary \ref{c:precomp} is also separable. Therefore, Theorem \ref{th:2main} also holds in that case. 

\section{\label{s:non-existence} Nonexistence of proper affine isometric actions}

It is easy to see that not every topological group $G$ admits a free continuous action by affine isometries on a Banach space $E$. For let $\pi$ be such an action, and let $\xi\in E$ be an arbitrary vector. The formula
\begin{equation}
d(g,h) = \norm{\pi_g(\xi)-\pi_h(\xi)}_E
\end{equation}
determines a left-invariant continuous metric on $G$. As a consequence, every point of $G$ has type $G_\delta$ (is an intersection of countably many open sets). Not every topological group has this property, and the easiest example would be a direct product, with Tychonoff product topology, of uncountably many copies of any non-trivial topological group, such as $\Z_2$ or the circle group $\T$. (Note that the group $\T^{\mathfrak c}$ is even separable, by a famous theorem of Pondiczery \cite{pondiczery} and Hewitt \cite{hewitt}.)

What about metrizable groups? Even in this case, the answer is in the negative. In fact, a typical sort of behaviour for a Polish group $G$ is to have non-compact stabilizers for actions by isometries on any {\em complete metric space.} 

Here we are going to explain this phenomenon.
Given a topological group $G$, denote by $\hat G^L$ the completion of $G$ with regard to the left uniform structure (Eq. (\ref{eq:left})). If $G$ is metrizable (which is the case we are interested in), then $\hat G^L$ is just the completion of $G$ equipped with any compatible left-invariant metric. 

It was noted early on by Dieudonn\'e \cite{dieudonne} that, while being a topological semigroup with jointly continuous multiplication (cf. \cite{RD}, Prop. 10.12(a)), the left completion $\hat G^L$ need not be a topological group. Among Dieudonn\'e's examples was the familiar group of homeomorphisms $\Homeo_+[0,1]$. 

\begin{example}[Dieudonn\'e \cite{dieudonne}]
The {\it right} completion of $\Homeo_+[0,1]$ is 
the semigroup of all continuous
order-preserving surjections from $[0,1]$ to itself (see also \cite{RD}, Exercise 1, p. 191), and so the left completion is formed by all relations on $[0,1]$ whose inverse relations are such surjections.
\end{example} 

\begin{example} 
The left completion of $U(\ell^2)$ is the semigroup of all linear isometric embeddings $\ell^2\hookrightarrow \ell^2$ (not necessarily onto) with the strong topology. 
\end{example}

\begin{example}
The infinite symmetric group $S_\infty$ consists of all self-bijections of a countably infinite set $\omega$ and is equipped with the Polish topology of pointwise convergence on $\omega$ with a discrete topology. The left completion of $S_\infty$ is the semigroup of all injections $\omega\hookrightarrow\omega$, with the pointwise convergence topology. 
\end{example}

\begin{example}
The left completion of the group $\Iso(\Ur)$ of isometries of the Urysohn space is the semigroup of all isometric embeddings of $\Ur$ into itself, with the pointwise topology. 
\end{example}

This list can go on and on.
Only in some special cases (e.g. where $G$ is locally compact, or abelian, or a Banach-Lie group) the one-sided completion of $G$ is a topological group again. However, the above sort of behaviour should not be regarded as ``pathological''. For instance, it serves as a basis for studying the famous distortion property / oscillation stability of the space $\ell^2$, as well as a number of other infnite Ramsey-type properties, in the context of topological transformation groups \cite{KPT}. (See also Chapter 8 in \cite{P06}, as well as the recent solution of the distortion problem for Urysohn sphere \cite{NVTS}.)

For a metric space $X$, denote by $\Emb(X)$ the semigroup consisting of all isometric embeddings of $X$ into itself (possibly proper ones), equipped with composition of maps as the binary operation and the topology of simple convergence on $X$ (induced by the embedding $\Emb(X)\hookrightarrow X^X$). A standard basic neighbourhood of an element $\tau\in\Emb(X)$ is of the form
\[V=V[\tau;\e;x_1,\ldots,x_n]=\{\tau^\prime\in\Emb(X)\colon d_X(\tau(x_i),\tau^\prime(x_i))<\e,~i=1,2,\ldots,n\}.\]

\begin{lemma}
The semigroup $\Emb(X)$ with the topology of simple convergence is a topological semigroup (with a jointly continuous multiplication). If the metric space $X$ is complete, then the semigroup $\Emb(X)$ is complete with regard to the left uniformity defined by the entourages of the diagonal
\begin{equation}
\label{eq:u}
V[\e;x_1,\ldots,x_n]_L=\{(\tau,\sigma)\in\Emb(X)^2\colon \forall i,~d_X(\tau(x_i),\sigma(x_i))<\e\}.\end{equation}
\end{lemma}

\begin{proof}
The first claim follows from the observation that, given
$\sigma,\tau\in\Emb(X)$, $\e>0$ and $x_1,x_2,\ldots,x_n\in X$, one has
\[V\left[\sigma;\frac{\e}2;\tau(x_1),\ldots,\tau(x_n)\right]\circ 
V\left[\tau;\frac{\e}2;x_1,x_2,\ldots,x_n\right]\subseteq
V[\sigma\tau;\e;x_1,x_2,\ldots,x_n].\]
The verification of the fact that the entourages in Eq. (\ref{eq:u}) form a basis for a uniform structure is pretty straightforward.
Now let $(\sigma_\alpha)$ be a Cauchy net with regard to the left uniformity. As an index, one can use a collection of the pairs of the form $(\e;F)$, where $\e>0$ and $F\subseteq X$ finite, partially ordered in an obvious way.
For every $x\in X$ the net $(\sigma_\alpha(x))$ is Cauchy in the metric space $X$ and so converges to a limit, which we denote $\sigma(x)$. For every $x,y\in X$ one has $d(\sigma_\alpha(x),\sigma_\alpha(y))=d(x,y)$, and the triangle inequality assures that $d(\sigma(x),\sigma(y))=d(x,y)$, that is, $\sigma$ as a map from $X$ to itself is an isometry and so $\sigma\in\Emb(X)$. It remains to verify that $\sigma_{\alpha}\to\sigma$, which is pretty much obvious.
\end{proof}

Notice that the restriction of left uniformity from the semigroup $\Emb(X)$ to the subgroup $\Iso(X)$ coincides with the left uniformity of the latter group. The left uniformity on $\Emb(X)$ is the only compatible uniformity admitting a basis consisting of entourages invariant under left translations. Defining in a similar way the right uniformity on $\Emb(X)$ is impossible. 

\begin{lemma}
The tautological action $\Emb(X)\times X\to X$ is jointly continuous.
\label{l:joint}
\end{lemma}

\begin{proof}
Let $\sigma\in\Emb(X)$ and $x\in X$, and let $\e>0$. If now $\tau$ and $y$ are such that $d(\sigma(x),\tau(x))<\e/2$ and $d(x,y)<\e/2$, then
\begin{eqnarray*}
d(\tau(y),\sigma(x))&\leq& d(\tau(y),\tau(x))+d(\tau(x),\sigma(x)) \\
&=& d(y,x) + d(\tau(x),\sigma(x)) \\
&<& \frac{\e}2+\frac{\e}2 =\e.
\end{eqnarray*}
\end{proof}

\begin{lemma}
Let a topological group $G$ act continuously by isometries on a metric space $X$. This action extends to a unique continuous action of the left completion $\hat G^L$ on $X$ by isometric embeddings.
\label{l:extends}
\end{lemma}

\begin{proof} Consider a continuous homomorphism $f\colon G\to\Iso(X)\hookrightarrow \Emb(X)$ determined by the action. The map $f$ is uniformly continuous with regard to the left uniformities, and since $\Emb(X)$ is a complete uniform space, $f$ extends to a unique uniformly continuous map $\hat f\colon\hat G^L\to\Emb(X)$. The continuity of the resulting map $G\times X\to X$ follows from Lemma \ref{l:joint}. It remains to verify that $\hat f(gh)(x)=\hat f(h)\hat f(g) (x)$. Choose left Cauchy nets $(g_\alpha)$, $(h_\alpha)$ in $G$ converging to $g$ and $h$, respectively. Then $g_\alpha h_\alpha\to gh$ because the multiplication in $\hat G^L$ is jointly continuous, and as $\hat f(g_\alpha h_\alpha)=\hat f(g_\alpha)\hat f(h_\alpha)$, the result follows.
\end{proof}

An element $z$ of a semigroup $S$ is {\em right cancellative} with regard to a subset $A\subseteq S$ if, whenever $x,y\in A$ and $xz=yz$, one has $x=y$.
An action of a topological group $G$ on a space $X$ is ({\em topologically}) {\em proper} if the inverses of compact sets under the map $G\times X\to X$ are compact. In particular, stabilizers of all points under a proper action are compact.

\begin{theorem}
\label{th:noprop}
Let $G$ be a topological group such that not every element of the semigroup $\hat G^L$ is right cancellative with regard to $G$. Then $G$ admits no free continuous action by isometries on a complete metric space. If moreover the stabilizer subgroup in $G$ of some element $z\in\hat G^L$ (under left multiplication) is non-compact, then every continuous action of $G$ by isometries on a complete metric space has non-compact stabilizers and in particular is not proper.
\end{theorem}

\begin{proof} Let $G$ act continuously by isometries on a complete metric space $X$. This action extends to a continuous action of the semigroup $\hat G^L$ on $X$ by isometric embeddings by Lemma \ref{l:extends} (cf. also Prop. 8.2.6 in \cite{P06}). Choose $x,y\in G$ and $z\in\hat G^L$ with the properties $x\neq y$ and $xz=yz$. 
One has for an arbitrary $\xi\in X$:
\[x(z(\xi))=(xz)(\xi)=(yz)(\xi)=y(z(\xi)),\]
meaning the action of $G$ on $X$ is not free. 
Under the additional assumption that the stabilizer of $z$ in $G$ is non-compact, we conclude that the stabilizer of $z(\xi)$ in $G$ is non-compact either, and so the action is not proper.
\end{proof} 

Numerous topological groups of importance have the properties listed in Theorem \ref{th:noprop}. 

\begin{theorem}
\label{th:concrete}
Each of the following Polish groups possesses a non-compact stabilizer of some element of the left completion, and therefore admits neither proper nor free continuous isometric actions on a complete metric space (in particular, neither proper nor free affine isometric actions on a Banach space):
\begin{itemize}
\item $U(\ell^2)$ with the strong operator topology,
\item $S_\infty$, with its standard Polish topology,
\item $\Iso(\Ur_1)$, with the topology of pointwise convergence,
\item $\Homeo_+[0,1]$, with the compact-open topology.
\end{itemize}
\end{theorem}

\begin{proof}
$G=U(\ell^2)$: let $v\in\hat G^L$ be a partial isometry whose image $\H=v(\ell^2)\subsetneq\ell^2$ has infinite codimension. The subgroup of $U(\ell^2)$ consisting of all block-diagonal operators of the form
\[\left(\begin{matrix} {\mathrm{Id}_{\H}} & 0 \\ 0 & u\end{matrix}\right),\]
where $u$ is a unitary operator on $\H^{\perp}$, stabilizes $v\in\widehat{U(\ell^2)}^L$ (when acting by multiplication on the left) and is isomorphic, as a topological group, to $U(\ell^2)$, hence non-compact.

For $S_\infty$, let $i$ be an injection $\omega\hookrightarrow\omega$ with infinite complement $\omega\setminus i(\omega)$, then $i$ is stabilized by a non-compact subgroup of $S_\infty$ consisting of all permutations whose restriction to $i(\omega)$ is identity.

For $\Iso(\Ur)$, it is enough to isolate two copies of $\Ur_1$ and embed them isometrically into $\Ur_1$ in such a way that they lie at a constant distance $1$ from each other and every self-isometry of the union extends to a global self-isometry of $\Ur_1$. Now an argument along the same lines applies. 

In the case of the group $G=\Homeo_+[0,1]$, it is more convenient to consider the action of $G$ {\em on the right} performed on the {\em right} completion $\hat G^R$ which consists, as we have seen, of all continuous order-preserving surjections of the interval to itself. Let $f$ be such a surjection which is not a homeomorphism. There is a $x\in [0,1]$ whose inverse image $f^{-1}(x)$ is a non-trivial interval. The {\em right} stabilizer of $f$ in $\Homeo_+[0,1]$ contains a subgroup isomorphic to $\Homeo_+[0,1]$, namely, a subgroup of self-homeomorphisms of the little interval $f^{-1}(x)$ preserving the endpoints and extended by the identity on the outside of $f^{-1}(x)$.
\end{proof}

The above observation is in sharp contrast with Haagerup and Przybyszewska's result \cite{HP} about locally compact groups admitting (metrically) proper affine isometric actions. (Recall that an action of $G$ on $X$ is {\em metrically proper} if, whenever $B$ is a bounded subset of $X$, the closure of the set $\{g\in G\colon gB\cap B\neq\emptyset\}$ is compact. This is an even stronger property than a topologically proper action.) 
It is one of those strange-looking yet typical properties of ``infinite-dimensional'' topological groups, such as extreme amenability \cite{P06,KPT,NVT}, Property $(OB)$ and its variations \cite{rosendal,atkin,hejcman,bergman}, ample generics and its various consequences \cite{KR}, etc.

In conclusion, we would like to correct a statement appearing in Gromov's book \cite{Gr}. Exercise (c) on p. 83 asks the reader to prove that every second countable group $G$ admits a free isometric action on the Urysohn space ($\Ur$, in our notation). As we have seen above (Theorem \ref{th:concrete}), many concrete Polish groups would not admit such an action. The exercise should ask the reader instead to prove that every second countable group $G$ admits an {\em effective} isometric action on the Urysohn space $\Ur$. This result belongs to Uspenskij \cite{Usp90,Usp98}.

\end{document}